 \documentclass[11pt]{article}
 \usepackage[top=1in,bottom=1in,left=1.5in,right=1.5in]{geometry}
  \usepackage{amsmath,amssymb}
  \usepackage{latexsym}
 \usepackage[dvips]{pstricks} 
  \usepackage{pst-node}
  \usepackage[dvips]{graphicx}

  \newcommand{\C}{\mathbb{C}}

  \newcommand{\e}{\mathbf{e}}
  \newcommand{\f}{\mathbf{f}}
  \newcommand{\g}{\mathbf{g}}
  \newcommand{\gl}{\mathbf{GL}}

  \renewcommand{\u}{\mathbf{u}}
  \renewcommand{\v}{\mathbf{v}}

  \newcommand{\x}{\mathbf{x}}
  \newcommand{\X}{\mathbf{X}}
  \newcommand{\y}{\mathbf{y}}
  \newcommand{\z}{\mathbf{z}}
  \newcommand{\0}{\mathbf{0}}

  \newcommand{\cX}{\mathcal{X}}
  \newcommand{\cY}{\mathcal{Y}}
  \newcommand{\cZ}{\mathcal{Z}}

  \newcommand{\hs}{\hspace*{\parindent}}
  \newcommand{\proof}{\hs \textbf{Proof.\ }}
  
  \newcommand{\tr}{\mathop{\mathrm{tr}}\nolimits}
  \newcommand{\adj}{\mathrm{adj\;}}

  \newcommand{\trans}{^\top}
  \newcommand{\qed}{\hspace*{\fill} $\Box$\\}

  \newcommand{\rank}{\mathrm{rank\;}}

  \newtheorem{theo}{\bfseries \hs Theorem}[section]

  \newtheorem{lemma}[theo]{\bfseries \hs Lemma}
  \newtheorem{claim}[theo]{\bfseries \hs Claim}
  \newtheorem{corol}[theo]{\bfseries \hs Corollary}

  \numberwithin{equation}{section} 

 \renewcommand{\span}{\mathrm{span}}

\begin{document}

 \title{A proof of the set-theoretic version\\of the salmon conjecture}

 \author
 {Shmuel Friedland and Elizabeth Gross\\
 Department of Mathematics, Statistics and Computer Science\\
 University of Illinois at Chicago\\
 Chicago, Illinois 60607-7045, USA\\
 \emph{email}: friedlan@uic.edu, egross7@uic.edu}

 \date{April 28, 2011}
 \maketitle

 \begin{abstract}
 We show that the irreducible variety of $4 \times 4 \times 4$ complex valued tensors of border rank
 at most $4$ is the zero set of polynomial equations of degree $5$ (the Strassen commutative
 conditions), of degree $6$ (the Landsberg-Manivel polynomials), and of degree $9$
 (the symmetrization conditions).

 \end{abstract}



 \section{Introduction}
 Let $V_r(m,n,l)\subseteq \C^{m}\otimes \C^n\otimes \C^l$ be the variety of tensors of border rank at most $r$. The projectivization of $V_r(m,n,l)$ is the $r$th secant variety of $\mathbb{P}^{m-1} \times \mathbb{P}^{n-1} \times \mathbb {P}^{l-1}$.  In 2007, Elizabeth Allman posed the problem of determining the ideal $I_4(4,4,4)$ generated by all polynomials vanishing on $V_4(4,4,4)$, this has been coined the salmon problem \cite{AR07}.
 The \emph{salmon conjecture}
 \cite[Conjecture 3.24]{PS} states that $I_4(4,4,4)$ is generated by polynomials of degree $5$ and $9$.
 A first nontrivial step in characterizing $V_{4}(4,4,4)$ is to characterize $V_4(3,3,4)$. In \cite{LM04}, Landsberg and Manivel show
 that $V_4(3,3,4)$ satisfies a specific set of polynomial equations of degree $6$.
 (See also \cite[Remark 5.7]{LM06} and \cite{OB10}.)
 Hence the revised version of the salmon conjecture
 states that $I_{4}(4,4,4)$ is generated by polynomials of degree $5$, $6$ and $9$ \cite[\S2]{Stu08}.
 This in particular implies the set-theoretic version of the salmon conjecture:
 $V_4(4,4,4)$ is the zero set of certain homogeneous polynomials of degree $5,6,9$.
 
 It is shown theoretically in \cite{Fr10} that $V_4(4,4,4)$  is cut out by polynomials of degree $5,9,16$.
 A main element toward this result is the characterization of the variety of tensors in $\C^{3}\otimes \C^3\otimes \C^4$ of border rank at most $4$. It is shown in \cite[Theorem 4.5]{Fr10} that $V_4(3,3,4)$ is cut out by polynomials of degrees
 $9$ and $16$ equations.  The degree $9$ equations follow from the observation in \cite{Fr10} that the four
 frontal slices of $\cX \in V_4(3,3,4)$, which consists of four $3 \times 3$ matrices, are symmetrizable by multiplication
 on the left and by multiplication on the right by a nonzero matrix $L,R\in \C^{3\times 3}$ respectively.
 The existence of nonzero matrices $L$ and $R$ is equivalent to the vanishing of all $9 \times 9$ minors of two corresponding$12 \times 9$ matrices whose entries are linear in the entries of $\cX$.  We call this set of polynomials the symmetrization conditions.

 Meanwhile, the entries of $L$ and $R$ are polynomials of degree $8$ in the entries of $\cX$.
 The degree $16$ equations in \cite{Fr10} are a result of the condition
 \begin{equation}\label{tracecond}
 LR\trans =R\trans L= \frac{\tr (LR\trans)}{3}I_3.
 \end{equation}
 The degree $16$ equations are used only in the case A.I.3 of the proof of Theorem 4.5 of \cite{Fr10}.

 In \cite{LM04}, Landsberg and Manivel exhibited ten linearly independent polynomials of degree $6$, referred here as the LM-polynomials, that vanish on $V_4(3,3,4)$.
 Using methods from numerical algebraic geometry, Bates and Oeding  give numerical confirmation that $V_4(3,3,4)$ is the zero set of specific set
 of polynomials of degree $6,9$ \cite{OB10}, where the polynomials of degree $6$ are the polynomials in \cite{LM04}.

 The aim of this paper to show that $V_4(3,3,4)$ is cut out by polynomials of degree $6$ and $9$.
 This is done by showing that in Case A.I.3 of \cite[Proof of Theorem 4.5]{Fr10} the use of polynomials of degree $16$ can be eliminated by use of the LM-polynomials.
 More precisely we show that any $3 \times 3 \times 4$ tensor $\cX=[x_{i,j,k}]\in\C^{3\times 3\times 4}$ whose four frontal are of the form
 \begin{equation}\label{334specten}
 X_k=\left[\begin{array}{ccc} x_{1,1,k}&x_{1,2,k}&0\\x_{2,1,k}&x_{2,2,k}&0\\0&0&x_{3,3,k}\end{array}\right], \quad k=1,2,3,4,
 \end{equation}
 has a border rank four at most if and only the ten LM-polynomials vanish on $\cX$.

 It is not too hard to show that a tensor $\cX\in \C^{3\times 3\times 4}$ of the form \eqref
 {334specten} has border rank four at most if and only if either the four matrices
 $\left[\begin{array}{cc} x_{1,1,k}&x_{1,2,k}\\x_{2,1,k}&x_{2,2,k}\end{array}\right], k=1,2,3,4$
 are linearly dependent or $x_{3,3,k}=0$ for $k=1,2,3,4$.
 Note that the condition that the above four $2 \times 2$ matrices are linearly dependent is equivalent
 to the vanishing of the polynomial
 \begin{equation}\label{4degreepol}
 f(\cX)=\det \left[\begin{array}{cccc}
 x_{1,1,1}& x_{1,2,1} &x_{2,1,1}& x_{2,2,1}\\
 x_{1,1,2}& x_{1,2,2} &x_{2,1,2}& x_{2,2,2}\\
 x_{1,1,3}& x_{1,2,3} &x_{2,1,3}& x_{2,2,3}\\
 x_{1,1,4}& x_{1,2,4} &x_{2,1,4}& x_{2,2,4}
 \end{array}\right].
 \end{equation}

 It turns out that the restrictions of ten LM-polynomials to $\cX$ of the form  \eqref
 {334specten} are the polynomials
 \begin{equation}\label{restLMpol}
 x_{3,3,k} x_{3,3,l}f(\cX) \textrm{ for } 1\le k\le l\le 4.
 \end{equation}
 Hence $\cX$ has a border rank four at most if and only the ten LM-polynomial vanish on $\cX$.
 Combining this with the results in \cite{Fr10} we deduce the the set-theoretic version of the salmon conjecture.

 We summarize briefly the content of the paper.  In \S2 we restate the characterization of $V_4(3,3,4)$ given in \cite[Theorem 4.5]{Fr10}.
 In \S3 we show that the use of polynomials of degree $16$ in the proof of \cite[Theorem 4.5]{Fr10} can replaced by the use of
 LM-polynomials.  In \S4 we summarize briefly the characterization of $V_4(4,4,4)$ as the zero set of polynomials of degree $5,6,9$.

 \section{A characterization of $V_4(3,3,4)$}

 We now state \cite[Theorem 4.5]{Fr10} which characterizes $V_4(3,3,4)$. Let $\cX=[x_{i,j,k}]_{i=j=k}^{3,3,4}\in \C^{3\times 3 \times 4}$. The four frontal slices of $\cX$ are denoted as the matrices $X_k=[x_{i,j,k}]_{i=j=1}^3\in \C^{3\times 3}, k=1,2,3,4$.
 Assume that $\cX\in V_4(3,3,4)$.
 A special case of \cite[Lemma 4.3]{Fr10}
 claims that there exist nontrivial matrices $L,R\in \C^{3\times 3}\setminus\{0\}$ satisfying the conditions
 \begin{eqnarray}\label{sufconsym1}
 L X_k-X_k\trans L\trans=0,\;k=1,\ldots,4, L\in\C^{3 \times 3},\\
 X_k R-R\trans X_k\trans =0,\;k=1,\ldots,4, R\in\C^{3\times 3}\label{sufconsym2}.
 \end{eqnarray}
These are the symmetrization conditions.

 If the entries of $R$ and $L$ are viewed as the entries of two vectors with $9$ coordinates each, then the systems \eqref{sufconsym1}
 and \eqref{sufconsym2} are linear homogeneous equations with coefficient matrices $C_L(\cX),
 C_R(\cX) \in \C^{12\times 9}$ respectively.  (Observe that for any $A\in \C^{3\times 3}$ the matrix $A-A\trans$ is skew symmetric, which has, in general, $3$ free parameters.)   The entries of $C_L(\cX),
 C_R(\cX)$ are linear functions in the entries of $\cX$.
 For a generic $\cX\in V_4(3,3,4)$, $\rank C_L(\cX)=\rank C_R(\cX)=8$.
 Hence we can express the entries of $L$ and $R$ in terms of corresponding $8 \times 8$ minors of $C_L(\cX),
 C_R(\cX)$ respectively.  There is a finite number of ways to express $L$ and $R$ in this way,
 and some of these expression may be zero matrices \cite{Fr10}. Thus the entries of $L$ and $R$ are polynomials of degree $8$ in the entries of in the entries of $\cX$.   If $\rank C_L(\cX)<8$ then each of this expression of $L$ is a zero matrix,
 and similar statement holds for $R$.  Hence if  $\rank C_L(\cX)=\rank C_R(\cX)=8$ then for each expression of $L$ and $R$
 the condition \eqref{tracecond} hold.   The characterization of $V_4(3,3,4)$ is given by \cite[Theorem 4.5]{Fr10}.
  \begin{theo}\label{334charbr4}$\cX=[x_{i,j,k}]_{i=j=k=1}^{3,3,4}\in\C^{3\times 3\times 4}$
 has a border rank $4$ at most if and only the following conditions holds.
 \begin{enumerate}
 \item\label{334charbr4a}
 Let $X_{k}:=[x_{i,j,k}]_{i=j=1}^3\in\C^{3\times 3},k=1,\ldots,4$ be the four frontal slices of $\cX$.
 Then the ranks of $C_L(\cX),
 C_R(\cX)$ are less than $9$.  (Those are $9-th$ degree equations.)
 \item\label{334charbr4b}
 Let $R, L$ be solutions of (\ref{sufconsym1})
 and (\ref{sufconsym2}) respectively be given by $8 \times 8$ minors of $C_L(\cX),
 C_R(\cX)$.  Then (\ref{tracecond}) holds.
 (Those are $16-th$ degree equations.)
 \end{enumerate}
 \end{theo}

 The proof of Theorem \ref{334charbr4} in \cite{Fr10} consists of discussing a number of cases.  The $16$ degrees polynomial conditions \eqref{tracecond} are used only in the case A.I.3.  In the next section we show how to prove the theorem in the case A.I.3 using only the ten
 LM-polynomials of degree $6$.

 \section{The case A.I.3 of \cite[Theorem 4.5]{Fr10}}

 Suppose $\cX\in \C^{3\times 3\times 4}$ and there exist two nonzero matrices $L,R\in \C^{3\times 3}$ such that (\ref{sufconsym1})--(\ref{sufconsym2})  hold.
 The case A.I.3 assumes that $L$ and $R$ are rank one matrices.  The degree $16$ equations yield that $LR\trans =R\trans L=0$, thus, the remainder of the proof of \cite[Theorem 4.5]{Fr10} in the case A.I.3 resolves the case where $LR\trans =R\trans L=0$.
 Therefore, to eliminate the use of polynomial conditions of degree $16$ we need to show the following.
 \begin{claim}\label{mclaim}  Let  $\cX\in\C^{3\times 3\times 4}$.  Let $R,L\in\C^{3\times 3}$ be rank one matrices satisfying the
 conditions  (\ref{sufconsym1})--(\ref{sufconsym2}) respectively.  Suppose furthermore that either $LR\trans\ne 0$ or $R\trans L\ne 0$.
 If the ten  LM-polynomials vanish on $\cX$ then $\cX\in V_4(3,3,4)$.
 \end{claim}

 In the rest of this section we prove Claim \ref{mclaim}.
 Assume that $L=\u\v\trans, R=\x\y\trans$.  The following claim is straightforward.
 \begin{eqnarray}\label{Lcond}
 \u\v\trans A \textrm{ is symmetric if and only if } \v\trans A= b\u\trans \textrm{ for some } b\in\C,\\
 A\x\y\trans \textrm{ is symmetric if and only if } A\x=c\y \textrm{ for some } c\in\C.  \label{Rcond}
 \end{eqnarray}
 By changing bases in two copies of $\C^3$ we can assume that $\u=\v=\e_3=(0,0,1)\trans$.
 Let $P,Q\in \gl(3,\C)$ such that
 \begin{equation}\label{PQconditions}
 P\trans\e_3,Q\trans\e_3\in \span(\e_3).
 \end{equation}
 Then if $A \in \C^3 \times \C^3$ such that $(\ref{Lcond})$ and $(\ref{Rcond})$ hold, $\e_3\e_3\trans (PA Q)$
  is symmetric.
 Observe next that $PAQ (Q^{-1}\x)(P\y)\trans$ is also symmetric.  Thus we need to analyze what kind of vectors can be obtained
 from two nonzero vectors $\x,\y$ by applying $Q^{-1}\x,P\y$, where $P,Q$ satisfy \eqref{PQconditions}.  By letting $Q_1:=Q^{-1}$
 we see that $Q_1$ satisfies the same conditions $Q$ in \eqref{PQconditions}.
 Hence $Q_1,P$ have the zero pattern
 \begin{equation}\label{zeropatPQ}
 \left[\begin{array}{ccc} *&*&*\\ *&*&*\\0&0&*\end{array}\right].
 \end{equation}

 \begin{lemma}\label{basnorm}  Let $\y\in\C^3\setminus\{\0\}$.
 If $\e_3\trans\y\ne 0$ then there exists $P\in\gl(3,\C)$ of the form \eqref{zeropatPQ} such that $P\y=\e_3$.
 If $\e_3\trans\y= 0$ then there exists $P\in\gl(3,\C)$ of the form \eqref{zeropatPQ} such that $P\y=\e_2$.
 \end{lemma}
 \proof  Assume first that $\e_3\trans\y\ne 0$.  Let $\f=(f_1,0,f_3)\trans,\g=(0,g_2,g_3)\trans\in\C^3\setminus\{\0\}$
 such that $\f\trans\y=\g\trans\y=0$.  Then $f_1g_2\ne 0$.  Hence there exists $P\in\gl(3,\C)$ of the form \eqref{zeropatPQ},
 whose first and the second rows are $\f\trans,\g\trans$ respectively, such that $P\y=\e_3$.

 Suppose now that $\e_3\trans \y=0$.  Hence there exists $P=P_1\oplus [1], P_1\in\gl(2,\C)$ such that $P\y=\e_2$.  \qed

 \begin{corol}\label{4choice}  Let $A\in\C^{3\times 3}$ and assume that $LA$ and $AR$ are symmetric matrices for some rank one matrices
 $L,R\in \C^{3\times 3}$.  Then there exists $P,Q\in\gl(3,\C)$ such that by replacing $A,L,R$ by $A_1:=PAQ, L_1:=Q\trans LP^{-1}, R_1=Q^{-1}RP\trans$ we can assume $L_1=\e_3\e_3\trans$ and $R_1$ has one of the following $4$ forms
 \begin{equation}\label{4choiceR}
 \e_3\e_3\trans, \quad \e_3\e_2\trans,\quad \e_2\e_3\trans,\quad \e_2\e_2\trans.
 \end{equation}
 \end{corol}

 To prove Claim \ref{mclaim}
 we need to consider the first three choices of $R_1$ in \eqref{4choiceR}.   Note that by changing the first two indices in $\cX\in\C^{3\times 3\times 4}$ we need to consider only the first two choices of $R_1$ in \eqref{4choiceR}.

 \subsection{The case $L=R=\e_3\e_3\trans$}
 Let $X_1,X_2,X_3,X_4\in\C^{3\times 3}$ be the four frontal sections of $\cX=[x_{i,j,k}]\in\C^{3\times 3\times 4}$.
 Assume that (\ref{sufconsym1})--(\ref{sufconsym2}) hold.  Then each $X_k$ has the form of \eqref{334specten}.
 (This is the case discussed in \cite[(4.7)]{Fr10}.)

 Using Mathematica, we took the 10 LM-polynomials available in the ancillary material of \cite[deg\_6\_salmon.txt]{OB10} and let $x_{1,3,k}=0$, $x_{2,3,k}=0$, $x_{3,1,k}=0$, $x_{3,2,k}=0$ for $k=1,2,3,4$.  The resulting polynomials had 24 terms.  We then factored $f(\cX)$ from these restricted polynomials.  This symbolic computations shows that
 that the restriction of the ten LM-polynomials to $\cX$ satisfying (\ref{sufconsym1})--(\ref{sufconsym2}) are the polynomials given in \eqref{restLMpol}.  Therefore, by the result of Landsberg-Manivel \cite{LM04}, if $\cX\in V_4(3,3,4)$
 then all polynomials in \eqref{restLMpol} vanish on $\cX$.

 Vice versa, suppose that all polynomials in \eqref{restLMpol} vanish on $\cX$.
 Assume first that the polynomial $f(\cX)$ given by \eqref{4degreepol} vanishes in $\cX$.
 Let
 \begin{equation}\label{defYk}
 Y_k=\left[\begin{array}{cc} x_{1,1,k}&x_{1,2,k}\\x_{2,1,k}&x_{2,2,k}\end{array}\right], \quad k=1,2,3,4,
 \end{equation}
 be the projection of the four frontal sections of $\cX$ given by \eqref{334specten} on $\C^{2\times 2}$.
 Then $f(\cX)=0$ if and only if $Y_1,Y_2,Y_3,Y_4$ are linearly dependent.  Assuming the generic case that $X_1,X_2,X_3,X_4$
 are linearly independent we can choose a new basis in $\span(X_1,X_2,X_3,X_4)$ of the form
 \begin{equation}\label{defZk}
 Z_k=\left[\begin{array}{ccc} z_{1,1,k}&z_{1,2,k}&0\\z_{2,1,k}&z_{2,2,k}&0\\0&0&0\end{array}\right], \; k=1,2,3,
 \quad Z_4=\e_3\e_3\trans.
 \end{equation}
 The tensor $\cZ=[z_{i,j,k}]\in\C^{3\times 3}$, whose three frontal sections are $Z_1,Z_2,Z_3$ is essentially $2\times 2\times 3$ tensor.
 Hence its border rank is at most $3$.  (Since $Z_1,Z_2,Z_3$ are linearly independent the rank of $\cZ$ is $3$.) Since $\rank Z_4=1$ we deduce that the border rank of $\cX$ is $4$ at most.  Therefore any $\cX$ of the form \eqref{334specten} which satisfies $f(\cX)=0$ has border rank $4$ at most.

 Assume now that $f(\cX)\ne 0$.  Since the ten polynomials in \eqref{restLMpol} vanish on $\cX$ it follows that
 $x_{3,3,k}=0$ for $k=1,2,3,4$.  In this case $\cX$ is essentially a $2\times 2\times 4$ tensor.  Hence its rank is $4$ at most.

 \subsection{The case $L=\e_3\e_3\trans, R=\e_3\e_2\trans$}
 Let $X_1,X_2,X_3,X_4\in\C^{3\times 3}$ be the four frontal sections of $\cX=[x_{i,j,k}]\in\C^{3\times 3\times 4}$.
 Assume that (\ref{sufconsym1})--(\ref{sufconsym2}) hold.
 This means that our tensor $\cX=[x_{i,j,k}]\in\C^{3\times 3\times 4}$ has the following zero entries
 $x_{1,3,k}=x_{3,1,k}=x_{3,2,k}=x_{3,3,k}=0$ for $k=1,2,3,4$.  So our tensor is essentially a $2\times 3 \times 4$ and  hence its border rank is $4$ at most \cite{Fr10}.

 \section{The defining polynomials of $V_4(4,4,4)$}
 In this section we state for reader's convenience the defining equations of $V_4(4,4,4)$.  We briefly repeat the arguments in \cite{Fr10}
 by replacing the degree $16$ polynomial conditions with the degree $6$ polynomial conditions.  Let $\cX=[x_{i_1,i_2,i_3}]\in \C^{4\times 4 \times 4}$.  For each $l\in\{1,2,3\}$ we fix $i_l$ while we let $i_p,i_q=1,2,3,4$ where $\{p,q\}=\{1,2,3\}\setminus\{l\}$.  In this way we obtain
 four $l$-sections $X_{1,l},\ldots,X_{4,l}\in \C^{4\times 4}$.  (Note that $X_{k,3}=[x_{i,j,k}]_{i=j=1}^4, k=1,2,3,4$ are the four frontal sections
 of $\cX$.)  Denote by $\X_l=\span(X_{1,l},\ldots,X_{4,l})\subset \C^{4\times 4}$ the $l$-section subspace corresponding to $\cX$.
 For each $l\in\{1,2,3\}$ we define the following linear subspaces of polynomials of degrees $5,6,9$ respectively in the entries of $\cX$.
 The defining polynomials could be any basis in each these linear subspaces.

 We first describe the Strassen commutative conditions.  (These conditions where rediscovered independently in \cite{AR03}.)  Take $U_1,U_2,U_3\in\X_l$.  View $U_i=\sum_{j=1}^4 u_{j,i}X_{j,l}$ for $i=1,2,3$.
 So the entries of each $X_{j,l}$ are fixed scalars and $u_{j,i}, i=1,2,3, j=1,2,3,4$ are viewed as variables.
 Let $\adj U_2$ be the adjoint matrix of $U_2$.
 Then the Strassen commutative conditions are
 \[U_1(\adj U_2)U_3-U_3(\adj U_2) U_1=0\]
 Since the values of $u_{j,i}, i=1,2,3, j=1,2,3,4$ are arbitrary, we regroup the above condition for each entry as a polynomial in $u_{j,i}$.
 The coefficient of each monomial in $u_{j,i}$ variables  is a polynomial of degree $5$ in the entries of $\cX$ and must be equal to zero.
 The set of all such polynomial of degree $5$ span a linear subspace, and we can choose any basis in this subspace.

 The degree $6$ and $9$ polynomial conditions are obtained in a a slightly different way.  Let $P=[p_{ij}],Q=[q_{ij}]\in \C^{4\times 4}$
  be matrices with entries viewed as variables.  View $PX_{k,l}Q,k=1,2,3,4$ as the four frontal section of $4\times 4\times 4$ tensor $\cX(P,Q,l)=[x_{i,j,k}(P,Q,l)]_{i,j,k=1}^4$.

  Let $\cY=[x_{i,j,k}(P,Q,l)]_{i,j,k=1}^{3,3,4}$. Now $\cY$ must satisfy the $6$ degree polynomial conditions of Landsberg-Manivel and the $9$ degree symmetrization conditions.  Since the entries of $P,Q$ are variables, this means that
 the coefficients of the monomials in the variables $v_{ij}, w_{ij}, i,j=1,2,3,4$ must vanish identically.  This procedure gives rise
 to a finite number of polynomial conditions of degree $6$ and $9$ respectively.  Again choose a finite number of linear independent
 conditions of degree $6$ and $9$ respectively.

 The zero set of the above polynomials of degrees $5,6,9$ defines $V_4(4,4,4)$.

 \section*{Acknowledgements}
 
 We thank Joseph Landsberg and Luke Oeding for helpful discussions regarding this problem.

\end{document}